\documentclass[12pt,a4paper]{article}
\usepackage[T1]{fontenc}
\usepackage[utf8]{inputenc}
\usepackage{lmodern}
\usepackage{textcomp}
\usepackage[english]{babel}
\usepackage{booktabs}
\usepackage{graphicx}
\usepackage{abstract}
\usepackage{indentfirst, comment}
\usepackage[colorlinks=true,linkcolor=blue,citecolor=blue,urlcolor=blue]{hyperref}
\usepackage{authblk}
\usepackage{appendix}
\usepackage{algorithmic}
\usepackage{lscape}
\usepackage[ruled, vlined]{algorithm2e}
\usepackage{multirow}
\usepackage{wrapfig}
\usepackage{subfigure}
\usepackage{float}
\usepackage{subfloat}
\usepackage{natbib}
\usepackage{fullpage}
 \usepackage{pgfplots}
\usepackage{xcolor}
\usepackage{xfrac}
\usepackage{rotating}
\usepackage{comment}
\SetKwRepeat{Do}{do}{while}%
\usepackage{amsthm,amsfonts,amsmath,amstext,amssymb,amsopn}
\usepackage{comment}
\usetikzlibrary{calc}

\usepackage{tikz}
\usepackage{setspace}
\usepackage{pgf}

\usetikzlibrary{spy}
\usetikzlibrary{arrows}
\usetikzlibrary{patterns,decorations.pathreplacing}
\usetikzlibrary{intersections}
\usepackage[font=small,labelfont=bf,labelsep=period,tableposition=top]{caption}
\usepackage{natbib}
\pgfplotsset{compat=1.12}
\usepackage{xcolor}
\usepackage{wrapfig}
\usepackage{framed}
\usepackage{pgfplots}
\usepackage{graphicx}
\usepackage{eurosym}
\usepackage{pgf-pie}
\usepackage{soul}
\usepackage{ulem}
\usepackage{float}
\usepackage{subfloat}
\usepackage{comment}

\newcommand{\GA}{\textsc{GA-CPR}}

\title{Gathering avoiding centralized pedestrian advice framework: an application for Covid-19 outbreak restrictions}
\author[1]{Veronica Dal Sasso}
\author[2]{Valentina Morandi}
\affil[1]{Optrail, Rome, Italy}
\affil[2]{Free University of Bozen - Bolzano, Faculty of Science and Technology, Italy}

\begin{document}

\maketitle

\begin{abstract}Due to the COVID-19 pandemic, the focus on everydays mobility has been shifted from traditional means of transport to how to safely commute for work and/or move around the neighbourhood. Maintaining the safe distance among pedestrian becomes crucial in big pedestrian networks. Looking at personal goals, such as walking through the shortest path, could lead to congestion phenomena on both roads and crossroads violating the imposed regulations. We suggest a centralized multi-objective approach able to assign alternative fair paths for users while maintaining the congestion level as lower as possible. Computational results show that, even considering paths that are not longer than the 1\% with respect to the shortest path for each pedestrian, the congestion phenomena are reduced of more than the 50\%.
\end{abstract}

\section{Introduction}

The steady progress towards a globalized world has, in the last decades, reduced the impact of distances: while in the past people were born, grew and spent all their lives in the same city, travelling has recently become more popular. The need to commute to reach the workplace 
consistently increased. As a consequence, the use of public and private means of transportation also increased, before the sudden drop due to the COVID-19 pandemic at the beginning of 2020. Before that, it was common to spend a lot of time travelling, ending up in traffic jams as every car driver tried to achieve his/her personal optimum, or being packed in overcrowded buses or trains. 
Now all this appears freezed, while the focus on everyday mobility is shifted on how to safely move around the neighborhood by foot. 
Walking also replaced public transports whenever possible, as people feel safer in the open air than in buses or metros. This means, however, that the choice of the path to walk through is seen in a new light, because keeping social distance is crucial also in the open air to reduce the risk of COVID-19 infections.
Clearly, looking at personal goals such as minimizing the length of your own path and being oblivious of others’ presence is self-defeating at keeping low rates of pedestrian congestion on the streets. On the other hand, by deferring the decisions to a centralized entity, which has a global view of the demand and capacity on the streets, all users could gain something in terms of experienced congestion phenomena.  

This idea was already mentioned by \cite{wardrop1952road}, where, speaking about vehicular traffic, the author makes the distinction between user equilibrium and system optimum: user equilibrium is a traffic assignment in which each user decides on its own the best route to follow while the system optimum is the traffic assignment in which the total travel time is minimized, without considerations on the behaviour and fairness among users. The deterioration of the overall solution in implementing the user equilibrium versus the system optimum is known in literature as the price of anarchy (see \cite{mahmassani1993network} and \cite{roughgarden2002bad} for details and mathematical background). A compromise solution between these two principles can offer interesting insights, leading to a win-win situation for all people involved. In fact, a miopic view of the problem from the perspective of the single user can lead to choices which, once put into practise, prove to be sub-optimal. As an example, let’s consider drivers commuting every day to and from the workplace. Each of them is prone to choose the shortest path, in terms of distance, travelling times or both. But, if their paths intersect, the level of congestion on the streets will increase the travelling times, reducing the effectiveness of the user’s choice. Even using real-time information on congestion, users' decisions simply result in a shift of congestion from previously congested roads to other roads. On the contrary, by choosing a different and at first sight less favourable path, the gain for the user may be substantial. As highlighted in \cite{avineri2009social},  however, the majority of users is not willing to act socially but instead selfishly, especially if the cost of acting socially is high \cite{fehr2003nature}. Different approaches for balancing user equilibrium and system optimum have been investigated  (see \cite{essen2016} and \cite{morandi2020bridging} for comprehensive reviews), spanning different ways of transportation. %
The bounded rational user equilibrium differs from the user equilibrium in the fact that the users are not completely free to choose their best route. In fact, a number of paths for each users can be considered according to the so-called indifference travel time band for which users would not feel the desire to change route. Details on the bounded rational user equilibrium can be found in \cite{mahmassani1987boundedly} and in \cite{zhang2011behavioral}. On the other hand, the constrained system optimum minimizes the total travel time trying to limit the unfairness among users by limiting the set of eligible paths only to those that are not longer than a small percentage with respect to the best choice. This centralized approach inspired this work. The first attempt to formulate the constrained system optimum, as a convex non-linear problem, can be found in \cite{Jahn2005}. For theoretical bounds on the price of anarchy we refer to \cite{Schulz2006}. Given the difficulty of handling non-linear latency functions, a first attempt to use a linear programming model to solve the constrained system optimum traffic assignment problem is proposed in \cite{angelelli2016proactive} and later in \cite{angelelli2020system}. Since the number of eligible paths is exponential, in \cite{dcvar2017} a fast and reliable heuristic algorithm to solve big road networks is proposed. In all presented approaches the set of eligible paths is generated a priori. Once the flow is routed on paths, it could happen that the experienced unfairness is much higher than the one evaluated a priori. To overcome the issue, in \cite{angelelli2020minimizing} two constrained system optimum formulations directly controlling the real experienced unfairness are presented. None of the presented approaches take into account the arc congestion level as a penalty for the system objective. The first attempt to embed arc congestion reduction techniques has been proposed in \cite{angelelli2019trade}, where a constrained system optimum with the aim of lowering congestion on worst congested arcs is proposed. 

A different approach which can be pursued to combine user equilibrium and system optimum is to formulate the problem as a multi-objective model. An example in the context of Air Traffic Management can be found in \cite{dalsasso2018}, where the authors assess the viability of incorporating single airlines preferences within a collaborative European framework. The main goal of the model presented is to ensure that air space capacity over Europe is never exceeded, while trying to accommodate stakeholders' requirements. Indeed, airlines request to the Air Traffic Manager trajectories and departure times for their flights, but usually they are not willing to share which policies lead to such choices. Hence, the requests may be disguised and may vary from one airline to the other, on the basis of the airlines' target. In order to adapt each airline's demand to the global objective of reducing costs, while ensuring the congestion level of air space sectors is below a set capacity, the model tries to minimize the deviation from the preferred trajectories, both in terms of delays and routing, and to minimize the total costs of traversing air sectors. The outcome of the model is a set of Pareto optimal solutions, among which one solution may be selected by further considerations on the fairness between the different stakeholders.

The problem of route assignment for pedestrians fits in these frameworks: the goal in this case is to ensure social distancing, and hence avoid congestion, while people walk on the streets. When the walking activity is motivated by the need of reaching a destination, most people would choose the shortest path available. Thus, quite intuitively, the users' satisfaction drops fast with the increase in the length of paths. Anyway, just a small increase may be necessary to avoid congestion. 

The aim of this paper is to show how paths may be assigned to pedestrians, balancing the opposing goals of minimizing the increase on the length for users and minimizing the capacity excess on both the streets and crossroads/squares. We recall that a path is considered eligible if it is no longer than a percentage of the shortest path, guaranteeing fairness among users. In Section \ref{compare}, the gathering avoiding pedestrian routing model is presented. In Section \ref{results}, the result of thorough experiments assessing the performance of the model on real road networks is shown. Finally, in Section \ref{conclusions}, some concluding remarks and future developments are provided.

\section{The gathering avoiding pedestrian routing model}\label{compare}

In this Section, we present the Gathering Avoiding Centralized Pedestrian Routing problem (GA-CPR). The model assigns paths to users in order to balance two opposing objectives, i.e. minimizing the increase in path length and minimizing the excess in the capacity. Hence, the problem is formulated as a multi-objective model. The fairness between users is embedded into the model by selecting as eligible paths only those that are not too unfair for users.

In the following, we will assume that the flow of pedestrians is constant. This allows us to neglect the time dimension of the problem and consider a static demand on the network. As assessed by \cite{sheffi1985urban}, this modelling choice holds when rush hour time slots are considered.

Let $G=\left(V,A\right)$ be a directed network, where $V$ represents the set of vertices and $A\subseteq V\times V$ the set of arcs. The set of arcs represents the set of roads in the pedestrian network while the set of vertices represents the set of crossroads between roads in set $A$. Each vertex $h \in V$ is associated to a capacity $cap_{h}$ representing the number of pedestrians that can transit through the vertex without causing gathering phenomena, and a traversing time $t_h$. Similarly, each road is associated with a capacity $cap_{ij}$, representing the number of pedestrians that can walk through that road segment without gatherings, and a traversing time $t_{ij}$. 
The requests for walking through the pedestrian area are collected and the demand of all pedestrians going from the same origin to the same destination are consolidated in an Origin-Destination (briefly OD). The set $C$ represents the set of OD pairs, each associated to an origin $O_c\in V$, a destination $D_c\in V$, and a pedestrian demand rate $d_c$ from $O_c$ to $D_c$. $t_{cSP}$ indicates the traversing time while walking along the shortest path between $O_c$ and $D_c$. Only paths that are similar to the shortest path in terms of length are given as input to the model for each OD pair. As it is necessary that the pedestrians feel that the proposed paths are among the best available. In details, we consider as eligible only paths whose relative excess in walking time with respect to the shortest path are within the fairness percentage $\phi$. The set of eligible pedestrian paths from origin to destination for each OD pair $c \in C$ is denoted by $K_c^{\phi}$. For each $k \in K_c^{\phi}$, let $t_{ck}$ be the time needed to go across path $k$. The indicator parameter $a_{ij}^{kc}$ has value 1 if path $k \in K_c^{\phi}$ traverses arc $\left(i,j\right) \in A$, while it has value 0 otherwise. Details on the generating paths algorithm are discussed and provided in \cite{angelelli2016proactive}.

Variables $x_{ij}$ represent the total pedestrian flow on arc $\left(i,j\right) \in A$, while variables $\sigma_{ij}$ represent the excess of flow with respect to arc capacity on arc $(i,j)$. 
Similarly, variables $\delta_{h}$ indicate the excess of flow traversing a certain vertex $h \in V$.  
Moreover, a number of variables are related to each path. Variables $y_{ck}$ represent the pedestrian flow of OD pair $c\in C$ routed on path $k \in K_c$.

The objective functions of the \GA\ model are denoted by
\begin{eqnarray}
    \tau(\phi) &=& \sum \limits_{c \in C}\sum \limits_{k \in K_c^{\phi}} \frac{t_{ck}}{t_{cSP}} y_{ck} \label{obj_time}\\
    \eta(\phi) &=& \sum \limits_{(i,j) \in A} \frac{t_{ij}}{cap_{ij}} \sigma_{ij} + \sum \limits_{h \in V} \frac{t_{h}}{cap_{h}} \delta_{h} \label{obj_excess}
\end{eqnarray}

Objective \ref{obj_time} records the total relative increase in walking time of pedestrians on paths with respect to the shortest one. It depends on the walking time unfairness parameter $\phi$. Objective \ref{obj_excess}, also depending on $\phi$, records the total relative excess of capacity for arcs and nodes weighted by the traversing time.

The \textsc{\GA} bi-objective model follows:
\begin{align}
\min \quad& \tau(\phi),\eta(\phi) \label{obj1}\\
& d_c = \sum \limits_{k \in K_c^{\phi}} y_{ck} & c \in C \label{MODEL:pathbased_1}\\
& x_{ij} = \sum \limits_{c \in C}\sum \limits_{k \in K_c^{\phi}} a_{ij}^{ck} y_{ck}
& (i,j)\in A \label{MODEL:pathbased_2}\\
& z_{h} = \sum \limits_{(i,j) \in A | j=h}x_{ij}
& h\in V \label{MODEL:pathbased_3}\\
& \sigma_{ij} \ge x_{ij}-cap_{ij}
& (i,j)\in A \label{MODEL:pathbased_4}\\
& \delta_{h} \ge z_{h}-cap_{h}
& h \in V \label{MODEL:pathbased_5}\\
& x_{ij} \ge 0 & (i,j) \in A \label{MODEL:pathbased_6}\\
& z_{h} \ge 0 & h \in V \label{MODEL:pathbased_7}\\
& \sigma_{ij} \ge 0 & (i,j) \in A \label{MODEL:pathbased_8}\\
& \delta_{h} \ge 0 & h \in V \label{MODEL:pathbased_9}\\
& y_{ck} \ge 0 & c \in C, k \in K_c. \label{MODEL:pathbased_10} \end{align}

Constraints (\ref{MODEL:pathbased_1}) ensure that the pedestrian demand $d_c$ of OD pair $c\in C$ is routed on paths in $K_c^{\phi}$. Constraints (\ref{MODEL:pathbased_2}) set the flow on an arc equal to the sum of the flows on the pedestrian paths traversing the arc. Constraints (\ref{MODEL:pathbased_3}) set the inflow of a vertex to be equal to the sum, among arcs entering in the vertex, of their pedestrian flows. Constraints (\ref{MODEL:pathbased_4}) set $\sigma_{ij}$, for each arc $(i,j) \in A$, greater or equal to the excess of flow $x_{ij}$ with respect to the arc capacity $cap_{ij}$. Notice that, because of constraints (\ref{MODEL:pathbased_8}) and objective function (\ref{obj1}), variable $\sigma_{ij}$ assumes value 0 whenever the arc capacity is not exceeded and $x_{ij}-cap_{ij}$, otherwise. Similarly, constraints (\ref{MODEL:pathbased_5}) set $\delta_{h}$ greater or equal to the excess of flow $z_{h}$ in vertex $h \in V$ with respect to the vertex  capacity $cap_{z}$. Together with Constraints (\ref{MODEL:pathbased_9}) and objective function (\ref{obj1}), variable $\delta_{h}$ assumes value 0 if the vertex capacity is not exceeded and $z_{h}-cap_{h}$, otherwise. Finally, constraints (\ref{MODEL:pathbased_6})-(\ref{MODEL:pathbased_10}) define the domain of the decision variables.

\subsection{Solution method}

A number of methodologies can be used to tackle multi-objective optimization. We will refer to \cite{chiandussi2012comparison} for a comprehensive review of such methodologies. 
We decided to use the linear combination of weights as solution method since, as assessed by \cite{chiandussi2012comparison}, the main advantages of this method are its simplicity and its efficiency (computationally speaking). However, its main drawback is the determination of the appropriate weight coefficients to be used in the final objective function. In fact, the choice of weights is crucial in determining the solution. Thus, the choice has to be made by the decision maker carefully considering the real-world problem characteristics. 
 
 The two objective functions $\tau(\phi)$ and $\eta(\phi)$ introduced in the previous paragraph are weighted according to the importance parameters $\alpha$ and 1-$\alpha$, respectively, with $\alpha \in [0,1]$.
 The resulting objective function is:
 
 $$\min \quad \alpha \sum \limits_{c \in C}\sum \limits_{k \in K_c^{\phi}} \frac{t_{ck}}{t_{cSP}} y_{ck} + (1-\alpha)[\sum \limits_{(i,j) \in A} \frac{t_{ij}}{cap_{ij}} \sigma_{ij} + \sum \limits_{h \in V} \frac{t_{h}}{cap_{h}} \delta_{h}]$$.

In Section \ref{results}, we will also analyze the impact of the choice of $\alpha$ parameter on the optimal solution. The aim is to give an overview of the impact of the two objective functions on different policy maker choices.

\section{Computational results}\label{results}

A benchmark of 4 map-based instances has been used in a computational study to assess the performance of the presented model. 
Instances are generated by randomly draw coordinates for origins and destinations of pedestrian flows from four Italian cities, namely Brescia, Bolzano, Rome and Vicenza. Arc walking times are obtained using Graphhopper and arc capacities are obtained by dividing the real distance on map by the safety walking distance of 2m imposed by Covid-19 regulations. Vertices' capacity is obtained as a percentage of the entering arcs capacity and vertices' traversing time is obtained randomly within a short time windows. OD pairs' demands are generated as a percentage of the capacity of the arcs exiting the origin. 
The 4 tested networks have 50 nodes with approximately 2500 arcs and 25 OD pairs, each one with a different demand. 
All instances can be found at: \url{https://valentinamorandi.it/research-outcomes/}.
For each instance, a traffic assignment has been found using a restricted path set with $\phi$ values in $\{0.01,0.05,0.1,0.15,0.2\}$ and $\alpha$ values in $\{1,0.9,0.7,0.5,0.3,0.1,0\}$. For each instance, we compute also the user equilibrium (in which each passenger goes on their shortest route) as a matter of comparison. It is obtained by setting $\alpha=1$ and $\phi=0$. In total, we obtain 36 traffic assignments for each instance.
The model is solved using CPLEX 12.6.0 on a Windows 64-bit computer with Intel Xeon processor E5-1650, 3.50 GHz, and 64 GB RAM.
 Results for the \GA\ model are presented and discussed in Sections \ref{models}. 

In the following all the computed and collected statistics are defined.

\begin{itemize}

\item \textbf{Congestion distribution}\\
\begin{itemize}
\item $\bar{\sigma}$: average relative excess of flow with respect to the arc capacity, i.e. $\frac{1}{|A|} \sum \limits_{(i,j) \in A} \frac{\sigma_{ij}}{cap_{ij}} $.
\item $\bar{\delta}$: average relative excess of flow with respect to the node capacity, i.e. $\frac{1}{|V|} \sum \limits_{h \in V} \frac{\delta_{h}}{cap_{h}} $.
\item $\lambda_{=0}$: Percentage of arcs and nodes with $\frac{\sigma_{ij}}{cap_{ij}}=0$ or $\frac{\delta_{h}}{cap_{h}}=0$ w.r.t. the total number of arcs and nodes.
\item $\lambda_{0<...< 0.25}$: Percentage of arcs and nodes with $0<\frac{\sigma_{ij}}{cap_{ij}}< 0.25$ or $0<\frac{\delta_{h}}{cap_{h}}< 0.25$ w.r.t. the total number of arcs and nodes.
\item $\lambda_{\ge 0.25}$: Percentage of arcs and nodes with $\frac{\sigma_{ij}}{cap_{ij}}\ge 0.25$ or $\frac{\delta_{h}}{cap_{h}}\ge 0.25$ w.r.t. the total number of arcs and nodes.
\end{itemize}

\item \textbf{User experience for each OD pair $c \in C$}
\begin{itemize}
\item $U^{ck}=\frac{t_{ck}-t_{SP}}{t_{SP}}$: user walking unfairness with respect to fastest path for the OD pair.
\item $\bar{U}=\frac{1}{\sum \limits_{c \in C} d_c} \sum \limits_{c \in C} \sum \limits_{k \in K_c^{\phi}} y_{ck} U^{ck}$: weighted user walking unfairness.
\end{itemize}

\item \textbf{Network statistics}
\begin{itemize}
\item $T$: increase of total walking time with respect to UE. 
\item $\Sigma$: reduction of the total walking time with exceeded capacity on arcs with respect to UE.
\item $\Delta$: reduction of the total walking time with exceeded capacity on nodes with respect to UE.

\end{itemize}
\end{itemize}

\subsection{Performance of the model}\label{models}

In this Section, we summarize the results obtained. Multi-objective models are interesting as they are able to show how contrasting goals are affected by each other. By increasing parameters $\phi$ and $\alpha$ we are, respectively, assessing the impact of loosening the constraint on users' fairness and of focusing on the reduction of congestion. 

\begin{figure}[htp]
\begin{tikzpicture}[scale=0.6]
\begin{axis}[width=7cm,height=6cm,ylabel={$\Delta$ \%},xlabel={$\phi$ (\%)},
ymin=-60,ymax=10,xmin=-2,xmax=22,
ymajorgrids=true,xtick pos=left,
title style={align=center}]
\addplot [solid, very thick] coordinates {
(0,0)
(1,-49.14)
(5,-49.23)
(10,-50.37)
(15,-51.5)
(20,-51.75)
};
\addplot [dotted, very thick] coordinates {
(0,0)
(1,-49.14)
(5,-49.23)
(10,-50.5)
(15,-51.61)
(20,-51.85)
};
\addplot [dashed, very thick] coordinates {
(0,0)
(1,-49.14)
(5,-49.17)
(10,-50.23)
(15,-51.33)
(20,-51.67)

};
\addplot [dashdotted, very thick] coordinates {
(0,0)
(1,-49.14)
(5,-49.17)
(10,-49.87)
(15,-51.13)
(20,-51.61)

};
\addplot [orange, very thick] coordinates {
(0,0)
(1,-49.14)
(5,-49.17)
(10,-49.85)
(15,-51.09)
(20,-51.6)

};
\addplot [orange, dotted,very thick] coordinates {
(0,0)
(1,-49.16)
(5,-49.19)
(10,-49.8)
(15,-50.22)
(20,-50.22)

};
\addplot [orange, dashed,very thick] coordinates {
(0,0)
(1,-2.24)
(5,-1.56)
(10,-2.31)
(15,-2.15)
(20,-3.75)

};

\end{axis}
\end{tikzpicture}
\quad\quad
\begin{tikzpicture}[scale=0.6]
\begin{axis}[width=7cm,height=6cm,ylabel={$\Sigma$ (\%)},xlabel={$\phi$ (\%)},legend pos=outer north east,
ymin=-70,ymax=10,xmin=-2,xmax=22, ymajorgrids=true,xtick pos=left,
title style={align=center}]
\addplot [solid, very thick] coordinates {
(0,0)
(1,-21.67)
(5,-31.39)
(10,-45.12)
(15,-55.42)
(20,-62.87)

};
\addlegendentry{$\alpha=0$}
\addplot [dotted, very thick] coordinates {
(0,0)
(1,-21.5)
(5,-27.65)
(10,-38.33)
(15,-46)
(20,-52.43)

};
\addlegendentry{$\alpha=0.1$}
\addplot [dashed, very thick] coordinates {
(0,0)
(1,-21.47)
(5,-27.34)
(10,-37.99)
(15,-45.68)
(20,-52.4)

};
\addlegendentry{$\alpha=0.3$}
\addplot [dashdotted, very thick] coordinates {
(0,0)
(1,-21.38)
(5,-27.22)
(10,-37.59)
(15,-45.01)
(20,-51.67)

};
\addlegendentry{$\alpha=0.5$}
\addplot [orange, very thick] coordinates {
(0,0)
(1,-21.2)
(5,-26.63)
(10,-37.12)
(15,-44.84)
(20,-50.43)

};
\addlegendentry{$\alpha=0.7$}
\addplot [orange, dotted,very thick] coordinates {
(0,0)
(1,-21.1)
(5,-26.54)
(10,-29.91)
(15,-31.26)
(20,-31.28)

};
\addlegendentry{$\alpha=0.9$}
\addplot [orange, dashed,very thick] coordinates {
(0,0)
(1,-1.74)
(5,-0.71)
(10,-1.71)
(15,-2.91)
(20,-2.08)

};
\addlegendentry{$\alpha=1$}

\end{axis}
\end{tikzpicture}
\caption{Decrease of walking time with exceeded capacity on nodes $\Delta$ (on the left) and on arcs $\Sigma$ (on the right), as a function of  $\phi$.}
\label{fig:gamma}
\end{figure}
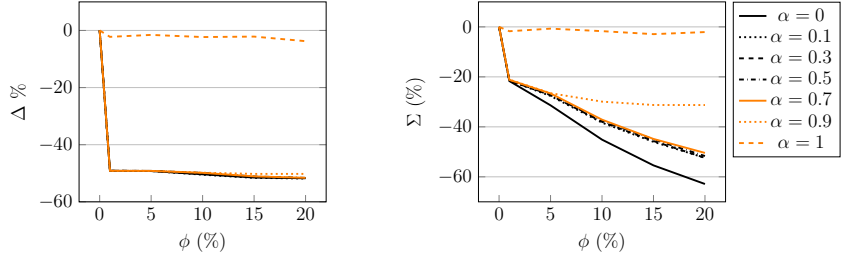

Quite intuitively, the walking time increases both when more importance is given to the reduction of congestion and when longer routes are made available (see Figure \ref{fig:tdrgtraveltime}). It is interesting to notice how, by just allowing a small freedom in the increase of the paths length, there is a considerable drop in the total time spent walking on congested nodes or arcs. Figure~\ref{fig:gamma} represents exactly this: the chart on the left shows how by, just adding those paths that are a 1\% longer than the shortest path, the reduction of walking time across congested nodes drops by 50\%. There is not much improvement, however, when longer paths are allowed. Another interesting results is that there is not much difference when the value of $\alpha$ changes (except, of course, when $\alpha=1$: when no weight is given to the congestion, there is no incentive of leaving the shortest path). On the right, a similar chart show an analogous results for the walking time on congested arcs. Notice that, in this case, the drop is more substantial when longer paths are allowed.

\begin{figure}[htp]
\begin{center}
\begin{tikzpicture}[scale=0.6]
\begin{axis}[ylabel={\%},xlabel={$\phi$ (\%)},
ymin=-1,ymax=8,xmin=-2,xmax=22,legend pos=outer north east,
ymajorgrids=true,xtick pos=left,
title style={align=center}]
\addplot [solid, very thick] coordinates {
(0,0)
(1,0.13)
(5,0.52)
(10,2.09)
(15,3.7)
(20,5.77)

};
\addlegendentry{$\alpha=0$}
\addplot [dotted, very thick] coordinates {
(0,0)
(1,0.13)
(5,0.33)
(10,1.28)
(15,2.41)
(20,3.82)

};
\addlegendentry{$\alpha=0.1$}
\addplot [dashed, very thick] coordinates {
(0,0)
(1,0.13)
(5,0.3)
(10,1.21)
(15,2.29)
(20,3.75)

};
\addlegendentry{$\alpha=0.3$}
\addplot [dashdotted, very thick] coordinates {
(0,0)
(1,0.13)
(5,0.3)
(10,1.15)
(15,2.25)
(20,3.65)

};
\addlegendentry{$\alpha=0.5$}
\addplot [orange, very thick] coordinates {
(0,0)
(1,0.13)
(5,0.29)
(10,1.13)
(15,2.24)
(20,3.3)

};
\addlegendentry{$\alpha=0.7$}
\addplot [orange, dotted,very thick] coordinates {
(0,0)
(1,0.13)
(5,0.29)
(10,0.52)
(15,0.66)
(20,0.66)

};
\addlegendentry{$\alpha=0.9$}
\addplot [orange, dashed,very thick] coordinates {
(0,0)
(1,0)
(5,0)
(10,0)
(15,0)
(20,0)

};
\addlegendentry{$\alpha=1$}
\end{axis}
\end{tikzpicture}
\end{center}
\caption{Increase in total walking time $T$, as a function of  $\phi$.}
\label{fig:tdrgtraveltime}
\end{figure}
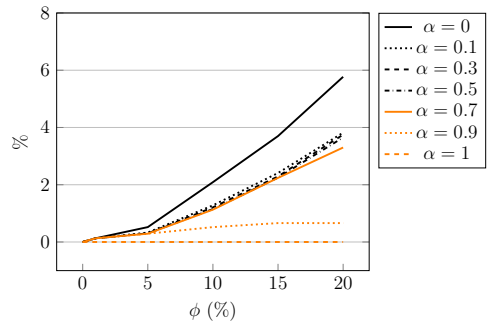

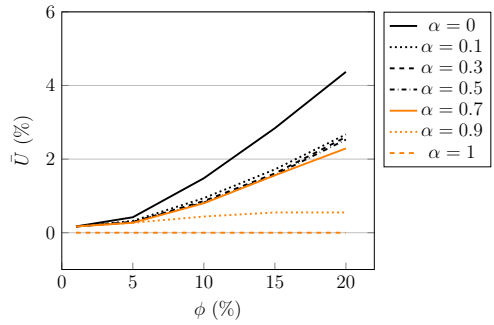
\begin{figure}[htp]
\begin{center}
\begin{tikzpicture}[scale=0.6]
\begin{axis}[ylabel={$\bar{U}$ (\%)},xlabel={$\phi$ (\%)},
ymin=-1,ymax=6,xmin=0,xmax=22,legend pos=outer north east,
ymajorgrids=true,xtick pos=left,
title style={align=center}]
\addplot [solid, very thick] coordinates {
(1,0.17)
(5,0.42)
(10,1.48)
(15,2.84)
(20,4.37)
};
\addlegendentry{$\alpha=0$}
\addplot [dotted, very thick] coordinates {
(1,0.17)
(5,0.32)
(10,0.94)
(15,1.72)
(20,2.67)
};
\addlegendentry{$\alpha=0.1$}
\addplot [dashed, very thick] coordinates {
(1,0.17)
(5,0.29)
(10,0.86)
(15,1.61)
(20,2.6)

};
\addlegendentry{$\alpha=0.3$}
\addplot [dashdotted, very thick] coordinates {
(1,0.17)
(5,0.29)
(10,0.82)
(15,1.57)
(20,2.53)

};
\addlegendentry{$\alpha=0.5$}
\addplot [orange, very thick] coordinates {
(1,0.17)
(5,0.27)
(10,0.8)
(15,1.56)
(20,2.29)

};
\addlegendentry{$\alpha=0.7$}
\addplot [orange, dotted,very thick] coordinates {
(1,0.17)
(5,0.27)
(10,0.44)
(15,0.55)
(20,0.55)

};
\addlegendentry{$\alpha=0.9$}
\addplot [orange, dashed,very thick] coordinates {
(1,0)
(5,0)
(10,0)
(15,0)
(20,0)
};
\addlegendentry{$\alpha=1$}

\end{axis}
\end{tikzpicture}
\end{center}
\caption{Average $\bar{U}$ user walking unfairness, as a function of  $\phi$.}
\label{fig:unf}
\end{figure}

Figure~\ref{fig:unf} shows how the average user unfairness, calculated as shown in the previous section, increases when longer paths are allowed. It clearly depicts the advantage of considering two contrasting objectives, as by varying the value of $\alpha$ the unfairness can be kept under control.

\begin{figure}[htp]
\begin{tikzpicture}[scale=0.6]
\begin{axis}[width=7cm,height=6cm,ylabel={$\bar{\delta}$ },xlabel={$\phi$ (\%)},
ymin=0,ymax=2,xmin=-2,xmax=22,
ymajorgrids=true,xtick pos=left,
title style={align=center}]
\addplot [solid, very thick] coordinates {
(0,1.65)
(1,0.76)
(5,0.75)
(10,0.72)
(15,0.7)
(20,0.69)

};
\addplot [dotted, very thick] coordinates {
(0,1.65)
(1,0.76)
(5,0.75)
(10,0.72)
(15,0.7)
(20,0.7)

};
\addplot [dashed, very thick] coordinates {
(0,1.65)
(1,0.76)
(5,0.76)
(10,0.73)
(15,0.71)
(20,0.7)

};
\addplot [dashdotted, very thick] coordinates {
(0,1.65)
(1,0.76)
(5,0.76)
(10,0.73)
(15,0.71)
(20,0.7)

};
\addplot [orange, very thick] coordinates {
(0,1.65)
(1,0.76)
(5,0.76)
(10,0.73)
(15,0.71)
(20,0.7)

};
\addplot [orange, dotted,very thick] coordinates {
(0,1.65)
(1,0.76)
(5,0.76)
(10,0.74)
(15,0.73)
(20,0.73)

};
\addplot [orange, dashed,very thick] coordinates {
(0,1.65)
(1,1.62)
(5,1.6)
(10,1.59)
(15,1.6)
(20,1.6)

};

\end{axis}
\end{tikzpicture}
\quad \quad
\begin{tikzpicture}[scale=0.6]
\begin{axis}[width=7cm,height=6cm,ylabel={$\bar{\sigma}$},xlabel={$\phi$ (\%)},
ymin=0,ymax=2.2,xmin=-2,xmax=22,legend pos=outer north east,
ymajorgrids=true,xtick pos=left,
title style={align=center}]
\addplot [solid, very thick] coordinates {
(0,1.99)
(1,0.14)
(5,0.17)
(10,0.19)
(15,0.26)
(20,0.16)

};
\addlegendentry{$\alpha=0$}
\addplot [dotted, very thick] coordinates {
(0,1.99)
(1,0.14)
(5,0.16)
(10,0.18)
(15,0.19)
(20,0.16)

};
\addlegendentry{$\alpha=0.1$}
\addplot [dashed, very thick] coordinates {
(0,1.99)
(1,0.14)
(5,0.16)
(10,0.18)
(15,0.18)
(20,0.15)

};
\addlegendentry{$\alpha=0.3$}
\addplot [dashdotted, very thick] coordinates {
(0,1.99)
(1,0.14)
(5,0.16)
(10,0.17)
(15,0.17)
(20,0.14)

};
\addlegendentry{$\alpha=0.5$}
\addplot [orange, very thick] coordinates {
(0,1.99)
(1,0.14)
(5,0.14)
(10,0.16)
(15,0.17)
(20,0.14)

};
\addlegendentry{$\alpha=0.7$}
\addplot [orange, dotted,very thick] coordinates {
(0,1.99)
(1,0.14)
(5,0.14)
(10,0.14)
(15,0.14)
(20,0.14)

};
\addlegendentry{$\alpha=0.9$}
\addplot [orange, dashed,very thick] coordinates {
(0,1.99)
(1,1.92)
(5,1.9)
(10,1.63)
(15,1.82)
(20,1.86)

};
\addlegendentry{$\alpha=1$}

\end{axis}
\end{tikzpicture}
\caption{Average node congestion $\bar{\delta}$ (on the left) and arc congestion $\bar{\sigma}$ (on the right), as a function of  $\phi$.}
\label{fig:sigmaavemax}
\end{figure}
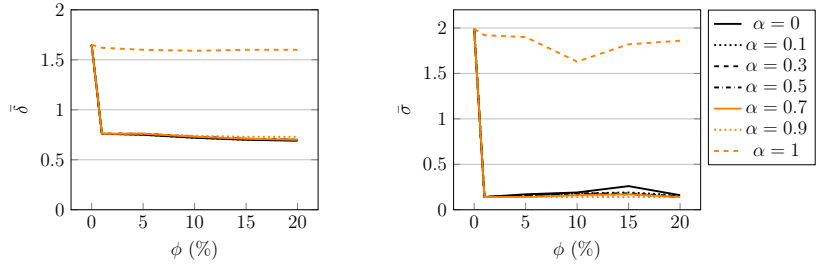

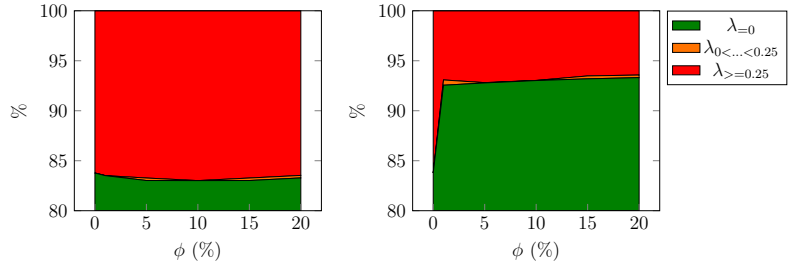
\begin{figure}
\begin{center}
\begin{tikzpicture}[scale=0.6]
\begin{axis}
[width=7cm,height=6cm,name=plot1,stack plots=y,
area style,ylabel={\%},xlabel={$\phi$ (\%)},
ymin=80,ymax=100,
xtick pos=left,
title style={align=center}]
\addplot [fill=green!50!black] coordinates
{(0,83.77)
(1,83.51)
(5,83.02)
(10,83.01)
(15,83.02)
(20,83.28)

}\closedcycle;
\addplot [fill=orange!90!red] coordinates
{(0,0)
(1,0)
(5,0.25)
(10,0)
(15,0.25)
(20,0.25)

}\closedcycle;
\addplot [fill=red] coordinates
{(0,16.23)
(1,16.49)
(5,16.73)
(10,16.99)
(15,16.73)
(20,16.47)

}\closedcycle;

\end{axis}
 \begin{axis}
[name=plot2,at={($(plot1.east)+(2cm,0)$)},anchor=west,height=6cm,width=7cm,name=plot2,stack plots=y,
area style,legend pos=outer north east,ylabel={\%},xlabel={$\phi$ (\%)},
ymin=80,ymax=100,
xtick pos=left,
title style={align=center}]
\addplot [fill=green!50!black] coordinates
{(0,83.77)
(1,92.55)
(5,92.78)
(10,93.02)
(15,93.19)
(20,93.31)

}\closedcycle;
\addlegendentry{$\lambda_{=0}$}
\addplot [fill=orange!90!red] coordinates
{(0,0)
(1,0.53)
(5,0.03)
(10,0.01)
(15,0.29)
(20,0.26)

}\closedcycle;
\addlegendentry{$\lambda_{0<...<0.25}$}
\addplot [fill=red] coordinates
{(0,16.23)
(1,6.92)
(5,7.19)
(10,6.97)
(15,6.52)
(20,6.43)

}\closedcycle;
\addlegendentry{$\lambda_{>=0.25}$}

\end{axis}
  
\end{tikzpicture}
\end{center}
\caption{Percentage of arcs and nodes in congestion classes with $\alpha=1$ and $\alpha=0$.}
\label{distrNode1}
\end{figure}

Figure~\ref{fig:sigmaavemax} remarks what already shown in Figure\ref{fig:tdrgtraveltime} but from another point of view. Instead of the total walking time on congested arcs and nodes, we are looking at the excess of flow. Congestion on nodes and arcs drops considerably when paths just 1\% longer than the shortest paths are taken into account. There is not much improvement, instead, in choosing longer paths, as expected. Figure~\ref{distrNode1} shows the congestion on arcs and nodes at the extreme points of the dominant solutions, i.e. the solutions of the multi -objective model obtained by setting $\alpha$ either $1$ or $0$. The percentage of non-congested arcs and nodes (in green) when only the total increase in walking time is minimized is very low, while it arises to almost 70\% when total increase in walking time is neglected and all the focus is on minimizing congestion. As noted also above, just a small increase on the length of the allowed paths gives a high reduction of congestion.

\section{Conclusions and future research} \label{conclusions}
In this paper we presented a multi-objective model for the Gathering Avoiding Centralized Pedestrian Routing Problem. Its aim is to ensure pedestrians can walk in safety according to COVID-19 regulations, which request people to social distancing. Hence, the objectives were to minimize the total increase in time needed by the users to reach their destinations and to minimize the exceeding of capacity on the arcs and at nodes. Taking into account users' fairness, the eligible paths are chosen among the ones that do not exceed the length of the shortest path of more than a proper percentage. Because of the contrasting nature of these objectives, the problem is modeled as a multi-objective model. 

The computational experiments were carried on by combining the two objectives linearly. The results are promising and show that just a small increase in the length of the paths that are taken into account can lead to a substantial reduction in congestion. Moreover, by opportunely setting parameter $\alpha$ it is possible to ensure the levels of unfairness among users is controlled. 

Future researches may improve the solution method by exploring different methods for multi-objective models. Heuristic methods to handle very large instances could be explored too.

\newpage
\bibliographystyle{natbib}
\bibliography{references}
\end{document}